\documentstyle[12pt]{article}
\def\thebibliography#1{\section*{\protect\large\bf References\markboth
 {REFERENCES}{REFERENCES}}\list
 {[\arabic{enumi}]}{\settowidth
\labelwidth{[#1]}\leftmargin \labelwidth
 \advance\leftmargin\labelsep
 \usecounter{enumi}}
 \def\newblock{\hskip .11em plus .33em minus -.07em}
 \sloppy
 \sfcode`\.=1000\relax}
\makeatletter
\@addtoreset{equation}{section}
\makeatother
\title{Structural properties of  weak cotype 2 spaces}
\author{Piotr Mankiewicz \and Nicole Tomczak-Jaegermann}
\newcommand\address{\noindent\leavevmode%
Institute  of Mathematics,\\
Polish Academy of Sciences\\
\'{S}niadeckich 8, \\
00-950 Warsaw, Poland,\\
{\small\tt% 
  piotr@impan.impan.gov.pl}\\[.5cm]
\noindent
Department of Mathematics,\\
University of Alberta,\\
Edmonton, Alberta, Canada T6G 2G1,\\
{\small\tt% 
  ntomczak@vega.math.ualberta.ca} }
\date{}

\newtheorem{thm}{Theorem}[section]
\newtheorem{prop}[thm]{Proposition}
\newtheorem{lemma}[thm]{Lemma}
\newtheorem{cor}[thm]{Corollary}

\newcommand{\rem}{{\bf Remark{\ \ }}}
\newcommand{\proof}{{\bf Proof{\ \ }}}
\newcommand{\qed}{\bigskip\hfill\(\Box\)}
\newcommand{\Rn}[1]{\mbox{{\it I\kern -0.25emR}$\sp{\,{#1}}$}}
\newcommand{\Nn}{\mbox{{\it I\kern-.16em\it N}}}
\newcommand{\di}{\hbox{\rm d}}
\newcommand{\ep}{{\varepsilon}}
\newcommand{\la}{{\lambda}}
\newcommand{\om}{{\omega}}
\newcommand{\Om}{{\Omega}}
\newcommand{\al}{{\alpha}}
\newcommand{\de}{{\delta}}
\newcommand{\be}{{\beta}}
\newcommand{\bc}{\mathop{\hbox{\rm bc}}}

\newcommand{\rank}{\mathop{\hbox{\rm rank}}}
\newcommand{\mix}{\mathop{\hbox{\rm Mix}}}
\newcommand{\Id}{\mathop{\hbox{\rm Id}}}
\newcommand{\vr}{\mathop{\hbox{\rm vr}}}
\newcommand{\codim}{\mathop{\hbox{\rm codim}}}

\newcommand{\vol}{\mathop{\hbox{\rm vol}}}

\newcommand{\Mix}{\mathop{\hbox{\rm Mix}}}

\newcommand{\cf}{{\it cf.}}
\newcommand{\cfeg}{{\it cf.\,e.g.,}}

\newcommand{\eg}{{\it e.g.,}}
\newcommand{\cfalso}{{\it cf.\,also}}
\newbox\nrmbox
\setbox\nrmbox=\hbox{$\Vert$}
\def\nrmrule{\vrule height\ht\nrmbox depth1.2\dp\nrmbox}
\setbox\nrmbox=%
  \hbox{\kern0.15em\nrmrule\kern0.15em\nrmrule\kern0.15em\nrmrule\kern0.15em}
\newcommand{\Tnorm}[1]%
  {\copy\nrmbox#1\copy\nrmbox}
\begin{document}
{\def\thempfn{}
\footnotetext{Supported in part by  KBN  (Poland),  NSERC (Canada)
and CRF (Canada, University of Alberta).}
\maketitle
\begin{abstract}
Several characterizations of weak cotype 2 and weak Hilbert spaces
are given in terms of basis constants and other structural
invariants of Banach spaces. For finite-dimensional spaces,
characterizations depending  on  subspaces  of fixed proportional 
dimension are proved. 
\end{abstract}
\section{Introduction}

Results of this paper concern weak cotype 2 spaces and weak Hilbert
spaces. Both classes are important in the local theory of Banach
spaces, by virtue of their connection to the existence of large
Euclidean  subspaces.
Recall that spaces of weak cotype 2 have been introduced by
V.~D.~Milman and G.~Pisier in \cite{MP} as spaces such
that every finite dimensional subspace contains a further subspace of
a fixed proportional dimension which is well Euclidean.
This class has  numerous other  characterizations, either by 
geometric  invariants  combined  with linear  structure,
or by inequalities between various ideal norms 
of related operators and their  $s$-numbers 
(see \eg\/ \cite{P} and references therein).

On the other hand, recent results of the authors (\cite{MT1},
\cite{MT3}) relate certain structural invariants of proportional 
dimensional  quotients of a finite-dimensional space
% , such as the Banach--Mazur distance or the basis constant, 
to volumetric invariants of the space; thus, via the well developed
theory, to the existence of Euclidean subspaces or quotients.
For example,  finite-dimensional results in \cite{MT3} imply,
although it is not explicitely stated in the paper,
that if a Banach space  has the property that its all subspaces
have a basis with a uniform upper bound for the basis constant, then
the space is of weak cotype 2. 
This property is obviously much too strong to 
characterize Banach spaces  of weak cotype 2;
for instance, spaces $L_p$ (with $p \ne 2$) contain subspaces 
without approximation property, hence without basis.
A natural question then arises whether a  weaker condition involving 
structural invariants of the same type
can in fact characterize spaces of
weak cotype 2. To put it more precisely, whether it is possible to
replace in the original weak cotype 2 definition, the property  of
being ``well Euclidean'' by   (much weaker)  properties of 
having some structural invariants  ``well bounded''.

One of the results of the present paper  (Theorem~\ref{rwaga} and the 
remark after Theorem~\ref{nwaga})  shows that this is indeed possible.
There exists  $\de_0 >0$ such that a Banach space $X$ is of weak cotype
2 if and only if every finite dimensional subspace $E$ contains a
a subspace $E_0 \subset E$ with $\dim E_0 \ge \de_0 \dim E$
such that   certain structural invariants (of $E_0$) have a
uniform upper bound. The  invariants  considered here
are the basis constants or the
complexification constants (for Banach spaces over reals only) or
the symmetry constants or in fact some other related parameters.

In the theory of proportional-dimensional subspaces of
finite-dimensional spaces is is sometimes of interest to deduce
properties of a fixed $n$-dimensional space $X$ from an information on
its all $\alpha n$-dimensional subspaces, with a fixed proportion 
$0 <  \alpha < 1$.  A fundamental example is well known and follows 
from the theory of type and cotype:  if 
all $\alpha n$-dimensional subspaces of $X$ are $C$-Euclidean than
$X$ itself  is $f(C)$-Euclidean (\cf\eg\/ \cite{T}). 
A recent more difficult  example can be found
in \cite{B} and  \cite{MT1}.
In Section 5 we study $n$-dimensional spaces $X$ such that the basis
constant of an arbitrary $\alpha n$-dimensional subspace $ E\subset X$
satisfies $\bc (E) \le C$.
This leads to proportional dimensional versions
of a result  mentioned above which follows  from \cite{MT3}.
In particular we show that the bound $\bc(E), \bc(F) \le C$ for
all $\alpha n$-dimensional subspaces $E$ and 
all $\alpha n$-dimensional
quotients $ F $ of $ X$ implies  that $X$ is a
weak Hilbert space, with the  constant bounded above by a
function of $C$.

The arguments in the paper are  based on two related ideas. The first
one is a technique developed  in \cite{MT1} and \cite{MT3} of finding
rather strange finite dimensional subspaces in spaces which fail to
have weak cotype 2. In the dual setting, which
is more convenient to use, it  can be described as follows.
First, for  a finite dimensional Banach space $E$,
using  deep facts from the local theory of Banach spaces,
we find a quotient $F$ which can be placed in a special 
position  in $\Rn{N}$;
and next we use a random (probabilistic) argument in order to prove
that majority of quotients of $F$ enjoys relative lack of well bounded
operators. The second idea comes from \cite{MT2}, where the 
authors have proved that a random 
proportional-dimensional quotient $F$ of $l_1^n$
cannot be embedded into 
a Banach space $F_1$  with a nice Schauder basis
and  $ \dim F_1 \le (1 + \de) \dim E$ for some small fixed $\de$.

The paper is organized as follows. In Section 2 we collect a
background  material related to  geometry and local theory of 
Banach spaces. In Section 3  we discuss spaces with  few 
well bounded operators and related volumetric lower  estimates.
The  main results of the paper  are  proved in Sections 4 and 5.
Section 6  contains a proof of a random result which generalizes
the result from \cite{MT2} to arbitrary finite-dimensional 
Banach spaces.

We shall consider only Banach spaces over reals. The complex case can
be delt with in an analogous manner. Our notation will follow  
\cite{P} and \cite{T}. We also refer the reader to 
\cite{P} for more details on  Banach spaces of  weak cotype 2.

\section{Preliminaries}

We will use the following geometric definition
of the weak cotype 2 spaces. A Banach space $X$
is of weak cotype 2 whenever there exist
$0 < \de_0 <1$ and $D_0 \ge 1$ such that
every finite-dimensional subspace $E $ of $X$
contains a subspace $\tilde{E} \subset E$
with $\dim  \tilde{E} = k = \ge \de_0 \dim E$
and with the Banach--Mazur distance satisfying
$\di (\tilde{E}, l_2^k) \le D_0$.
This definition is equivalent to the one most commonly used at present;
in fact, it is shown in \cite{P}, Theorem 10.2,
that the weak cotype 2 constant of $X$
satisfies
\begin{equation}
wC_2(X) \le C \de_0^{-1} D_0,
\label{weak_cot_2}
\end{equation}
where $C$ is a universal constant.

Recall that  a Banach space $X$ is of weak type 2
whenever $X^*$ is of weak cotype 2
and $X$ is $K$-convex. In particular, the
weak type 2 constant satisfies
\begin{equation}
wT_2(X) \le K(X) wC_2(X^*).
\label{weak_typ_2}
\end{equation}
A Banach space $X$ is a  weak Hilbert space
if $X$ is  of weak type 2 and  of weak cotype 2.
We have no use of the technical
definition of the weak Hilbert constant,
let us just recall that this constant is controlled
from above by $wT_2(X)$ and $wC_2(X)$.
Finally, the following inequality is an easy consequence of
the result of Pisier from \cite{PClev}, Theorem 2 and
Corollary 9, which in turn are related to Pisier's deep
$K$-convexity theorem.  If $X$ is a weak Hilbert, then the
$K$-convexity constant of $X$ satisfies, for any 
$0 < \theta \le 1$,
\begin{equation}
K(X) \le C(\theta) (wT_{2}(X) wC_{2}(X))^{\theta}.
\label{kconvex}
\end{equation}

Let $(E,\|\cdot\|)$ be an $n$-dimensional Banach space. 
For the Banach--Mazur distance $\di (E, L_2^n)$ from $E$
to the Euclidean space $l_2^n$ we will use a shorter notation
of  $\di_E$.
Fix a Euclidean norm $\|\cdot\|_2$ on $E$ and 
identify  $E$ with $\Rn{n}$ in such a way that $\|\cdot\|_2$
becomes the natural  $l_2$-norm on $\Rn{n}$.
Let us recall that the volume ratio of $E$, $\vr(E)$, is defined by
$$
\vr (E)  =  \bigl(\vol B_E / \vol {\cal E}_{\max} \bigr)^{1/n},
$$
where  ${\cal E}_{\max} \subset B_E $ is  the ellipsoid of
maximal volume contained in $B_E$.

More generally,  for any ellipsoid ${\cal E} \subset B_E $,
let $|\cdot|_2$ be the associated Euclidean norm, and let
$\rho = (\vol B_E / \vol {\cal E} )^{1/n}$.
Szarek's volume ratio result 
(\cfeg\ \cite{P} Theorem 6.1) says that 
for any $1 \le k < n $, there
exists a subspace ${H}  \subset E$ with ${\dim {H}}= k$
such that
\begin{equation}
c \rho^{- n/(n-k)} |x|_2 \le \|x\| \le |x|_2
\quad \mbox{for }x \in H,
\label{vol_rat}
\end{equation}
where $c >0 $ is a universal constant.
It should be mentioned that some other volumetric invariants
allow  estimates with  much better
asymptotic dependence  on $\lambda = k/n$, as $\lambda \to 1$
(\cfeg\ \cite{P}). However application
of these more delicate methods would complicate proofs without
making essential improvements to  final inequalities.

For $k \le n$ set
$$
V_k (E) = V_k (B_E)
=  \sup \{(\vol P_F(B_E) / \vol  P_F({B_2^n}))^{1/k}\, |\,
F \subset E, \dim F =k\}.
%\label{vol-inv}
$$
This invariant is related to the notion of volume numbers
of operators (\cf \cite{P} Chapter 9). In particular,
$V_k (E) \le V_l (E)$ for $1 \le l \le k \le n$.

The relevance of this invariant to the problem of finding
Euclidean quotients of $E$ is described by the following
standard lemma. Its proof is based on 
Santal\`{o} inequality and the volume ratio method
(\ref{vol_rat}) used in the dual space $E^*$.
We leave further details to the reader.

\begin{lemma}
\label{volumes}
Let $0 < \beta  <1$.
Let $E = (\Rn{n}, \|\cdot\|)$ be a Banach space
such that $B_E \subset B_2^n$ and
let $a >0$ satisfies  $V_{\beta n}(E) \ge a$. For every
$\sigma >0$ there is a quotient $G$ of $E$ such
that $\dim G \ge  \beta \sigma (1 + \sigma)^{-1} n $
and $\di_G \le C a^{-(1 + \sigma)}$, where
$C \ge 1$ is a universal constant.
\end{lemma}

The following fact is an easy consequence of \cite{P}, Lemma 8.8.
Let $\cal E$ be an ellipsoid on $E$ and $F$ be a quotient of $E$
with $\dim F = \lambda n$
and with the quotient map $Q: E \to F$.  Then 
\begin{equation}
\bigl( \vol Q(B_E)/ \vol Q({\cal E})\bigr)^{1/\lambda n}
     \le  a(\lambda) 
      \bigl( \vol B_E/ \vol {\cal E}\bigr)^{ 1/n},
  \label{prop_vol}
\end{equation}
where $a(\lambda) \ge 1$ depends on $\lambda$ only. In particular,
$\vr (F) \le a(\lambda) \vr (E)$.

The next lemma shows that given a finite-dimensional Banach space
one  can dramatically improve  geometric properties of its  unit ball
by passing to  quotients of proportional dimensions. 
The argument is based on
several  deep results in the local theory
of Banach spaces (\cite{Mi}, \cite{BS}, \cfalso\ \cite{P})..
The  lemma implicitly uses the ellipsoid
of minimal volume containing the unit ball $B_E$ of a given
$n$-dimensional space $E$; in fact it is concerned
with the following  property of a Euclidean norm   $|\cdot|_2$
on $E$:  there exists $c > 0$ such 
that every  rank $k$ orthogonal projection $P$ in $(E,|\cdot|_2)$
satisfies
\begin{equation}
\|P: E \to (E, |\cdot|_2)\| \ge c (k/n)^{1/2}.
\label{projec'}
\end{equation}
The norm $\Tnorm{\cdot}$ associated to the minimal volume
ellipsoid satisfies  (\ref{projec'})
with  $c=1$ (\cfeg\ \cite{T} Proposition 3.2.10).

Note that if $E$ is a finite-dimensinal Banach space and  
$|\cdot|_2$ is an Euclidean
norm on $E$ then each qoutient map $Q\,: E \to F$ induces on $F$ in a
natural way an Euclidean norm  $|\cdot|_{2,F}$. Namely, we set
$ B_{2,F} = Q(B_{2,E})$, where $B_{2,F}$ and $B_{2,F}$ stand for
Euclidean balls in $F$ and $E$ respectively.   

\begin{lemma}
\label{quotients}
For every $0 < \lambda <1$ there is $\rho = \rho (\lambda)\ge 1$
and for every $c>0$ there is  
$\tilde{\kappa} = \tilde{\kappa} (\lambda, c)\ge 1$,
such that  the following conditions hold for any 
$n$-dimensional Banach space   $E$:
\gdef\labelenumi{\theenumi}
\gdef\theenumi{(\roman{enumi})}
\begin{enumerate}
\item  
There exists a Euclidean norm $|\cdot|_2$ on ${E}$, 
with the unit ball $B_2$, such that 
\begin{equation}
(2^{1/2}\di_E)^{-1}B_2 \subset B_{{E}}  \subset  B_2
  \label{dist_incl}
\end{equation}
and it satisfies (\ref{projec'}) with $c = 2^{-1/2}$;
\item
there exists a $\lambda n$-dimensional quotient $F$ of $E$
satisfying ${\vr}({F}) \le \rho$;
\item
if $|\cdot|_2$ is a Euclidean norm on ${E}$ satisfying
(\ref{projec'}) for some $c >0$, then
there exists a $\lambda n$-dimensional quotient $\tilde{F}$ of $E$
and  an orthonormal basis $\{x_i\}$ in $(\tilde{F}, |\cdot|_2)$
such that
\begin{equation}
\label{kappa}
\max_i \|x_i\|_{\tilde{F}} \le \tilde{\kappa}.
\end{equation}
\end{enumerate}
\end{lemma}
\proof
The proof of $(ii)$ and $(iii)$ was already given in 
\cite{MT1}, Proposition 3.5. 
As for  $(i)$,        
the Euclidean norm described 
in  this  condition   combines properties
of the  norm $\Tnorm{\cdot}$ associated with
the ellipsoid of minimal volume and a  norm which determines
the Euclidean distance $\di_E$. Indeed, let 
$\Tnorm{\cdot}'$ be a norm  satisfying
$$
\Tnorm{x}' \le \|x\| \le  \di_E \Tnorm{x}'
\quad {\rm for} \quad x \in E.
$$
Then for $x \in E$ set
$$
|x|_2 = 2^{-1/2}(\Tnorm{x}^2+ \Tnorm{x}'^2)^{1/2}.
$$
Clearly, $|x|_2 \le \|x\| \le 2^{1/2} \di_E |x|_2$, for $x \in E$,
hence (\ref{dist_incl}) holds. 
To prove that $|\cdot|_2 $  satisfies (\ref{projec'}),
observe that for every $x \in E$ and every rank $k$
orthogonal projection $P$ in $(E, |\cdot|_2 ) $ one has
\begin{eqnarray*}
|Px|_2 &= &2^{-1/2} (\Tnorm{Px}^2 + \Tnorm{Px}'^2)^{1/2}  \\
 &\geq &  2^{-1/2} \Tnorm{Px} \geq 2^{-1/2} \Tnorm{P_1Px}
              = 2^{-1/2} \Tnorm{P_1x},
\end{eqnarray*}
where $P_1$ is the orthogonal projection in $(E,\Tnorm{\cdot})$ 
with $\ker P_1 = \ker P$ and next apply (\ref{projec'})
in the space $(E,\Tnorm{\cdot})$.
\qed

If $F$ is a finite-dimensional space and $B_2$ is a Euclidean
ball on $F$, then, for a fixed orthonormal basis $\{x_i\}$ in $F$, 
by $B_1$ we shall denote $abs\,conv\,\{ x_i\}$. 

\begin{cor}
  \label{spec_pos_cor}
For every $0 < \lambda <1$ there is $\rho = \rho (\lambda)\ge 1$
and  $ {\kappa}=  {\kappa} (\lambda) \le 1$,
such that an arbitrary $n$-dimensional Banach space   $E$
satisfies:
\gdef\labelenumi{\theenumi}
\gdef\theenumi{(\roman{enumi})}
\begin{enumerate}
\item  
there exist a $\lambda n$-dimensional quotient $F$ of $E$
and a Euclidean ball $B_2$ on $F$ such that 
${\vr}({F}) \le \rho$ and 
$$
   {\kappa} B_1 \subset B_F
\qquad \mbox{and}
\qquad  (2^{1/2}\di_E)^{-1} B_2 \subset B_F \subset  B_2.
$$
\item  
there exist a $\lambda n/2$-dimensional quotient $F$ of $E$
and a Euclidean ball $B_2$ on $F$ such that for some 
$ (2^{1/2}\di_E)^{-1} \le  {a} \le 1$
we have
$$
   {\kappa}  B_1 \subset  B_F,
\qquad   {a} B_2 \subset B_F \subset  B_2,
\qquad  \bigl(\vol B_F/ \vol (a B_2)\bigr)^{1/\lambda n} \le \rho.
$$
\end{enumerate}
\end{cor}
\proof
Condition $(i)$ follows directly from Lemma~\ref{quotients}
by chosing a Euclidean norm on $E$ satisfying
(\ref{dist_incl}) and then 
passing twice to quotient spaces and using (\ref{prop_vol}).
Notice that a Euclidean ball $B_2$  on  $E$ determines the natural
Euclidean ball on every  quotient $F$ on $E$, and if $B_2$ satisfies
(\ref{projec'})  and $\dim F$ is proportional to $\dim E$
then the ball on $F$  satisfies (\ref{projec'}) as well,
with the  constant  depending on the proportion.

To get $(ii)$, fix $0 < \lambda' < 1$  to be determined later.
First pass to a $\lambda'n$- dimensional quotient $F'$ of $E$
satisfying Lemma~\ref{quotients} $(i)$ and $(ii)$. Let
$\cal E$ be the maximal volume ellipsoid on $F'$;
there exists a quotient $F''$ of $F'$ with
$\dim F'' = \dim F'/2$  and the quotient map $Q: F' \to F''$
such that on $F''$ we have
$Q({\cal E}) = a Q( B_2)$, for some $a$.
Clearly, $a \le 1$ and $(i)$ implies that 
$ a \ge (2^{1/2}\di_E)^{-1} $.
Passing to a quotient $F$ of $F''$ with $\dim F = \lambda' \dim F''$,
and using Lemma~\ref{quotients} $(iii)$,
we get all required inclusions; the bound
for the ratio of volumes follows from (\ref{prop_vol}).
Given $1/2 < \lambda < 1 $ choose $\lambda'$ such
that $\dim F = \lambda n$, then all constants
involved will depend on  $\lambda$.
\qed

Volumetric techniques for finding Euclidean sections provide
Euclidean subspaces of small proportional dimensions.
The next lemma describes a method from \cite{MP}
of constructing 
Euclidean subspaces of large proportional dimensions in spaces
saturated with small Euclidean ones. The proof of the statement below
can be found in \cite{MT1} Theorem 4.2.

\begin{lemma}
\label{bighilbert}
Let $ 0 < \delta < \xi < 1$.
Let $Z$ be an $N$-dimensional space such that every
$\xi N$-dimensional subspace
$Z_1$ of $Z$ contains a subspace $H$ with
$\dim H \ge (\xi - \delta) N$ such that
$ \di_H \leq \ D$,  for some  $D \geq 1$.
Then for every $ 0 < \eta < 1 - \delta $
there exists a subspace
$\tilde{H}$ of $Z$ with 
$\dim \tilde{H}  \ge (1  - \delta  - \eta )N$
such that
$ \di_{\tilde{H}}  \leq c D$,
where $c = c(\xi, \delta, \eta)$.
\end{lemma}

Let $E$ be  an $n$-dimensional Banach space and let
$ \| \cdot \|_2$ be a Euclidean norm on $E$.
Recall that an operator $T: E \to E$ is
said to be $(k, \beta)$-mixing for $k,\beta \ge 0$, if
and only if
there is a subspace $F \subset E$ with $\dim F \geq k$ such 
that $|P_{F^{\perp}}Tx|_2 \geq \beta |x|_2$ for every $x \in F$, 
where $P_{F^{\perp}}$ denotes the the orthogonal projection onto
$F^{\perp}$. If this is the case then we write 
$T \in \Mix_n(k, \beta)$.  The fact whether a fixed
operator $T$ is in $\Mix_n(k,\beta)$ may depend on the
choice of the Euclidean structure on $E$. 
Clearly, if $k \ge \ell$, then
$\Mix_n(k, \beta) \subset \Mix_n(\ell, \beta)$.
The following proposition is a folklore one
(\cfeg\ \cite{S2} Lemma 3.4A, \cite{M3}).

\begin{prop}
\label{mix}
Let $E$ be an $n$-dimensional Banach space. For
an arbitrary Euclidean norm $|\,.\,|_2$ on $E$ one has
\gdef\labelenumi{\theenumi}
\gdef\theenumi{(\roman{enumi})}
\begin{enumerate}
\item 
for every  projection $P$ of rank $k \leq n/2$,
we have $ 2P \in \Mix_n(k,1) $;
\item 
if the basis constant $\bc (E) \leq M$ then 
for every $ k \in \Nn$ there is an
operator $T: E \to E$ with $\|T\| \leq 2M$ and
$ T \in \Mix_n(k, 1)$;
\item 
let $T: E \to E$, $k \in \Nn$ and $\be \ge 0$. Then
$ T \in \Mix_n(k,\beta)$ if and only if 
$T^* \in \Mix_n (k,\beta)$ 
(with respect to the dual Euclidean norm $|\cdot|_2^*$).
\end{enumerate}
\end{prop}

\section{Volumetric estimates}
\label{fin_est}

Technical result which this paper is based upon yields the existence,
for a given finite-dimensional Banach space, of a quotient space, say
$F$, of proportional dimension which admits  relatively  few well bounded
operators.  Moreover the same property is satisfied in any further
quotient $F_0$ of $F$ and in any space $\tilde{F}$ which admits $F$ as
its quotient, provided that the dimension of the new spaces is close
enough to the dimension of $F$.

In this section we shall work with an $N$-dimensional Banach space 
$E = (\Rn{N}, \|\cdot\|)$, on which we always consider
the (natural) Euclidean norm $\|\cdot\|_2$ and the associated
Euclidean ball $B_2^N$. If $F$ is a quotient of $E$, with the quotient
map $Q: E \to F$, then the natural Euclidean norm on $F$ has the 
unit ball $Q(B_2^N)$. Unless otherwise stated, 
these natural Euclidean norms are used for all geometric invariants.

First we discuss quotients which admit lower estimates for norms of
mixing operators. 
The main analytic  estimate is stated in the  following  theorem;
the   proof  involves a  random construction and it
is  postponed until Section~\ref{random}.

\begin{thm}
\label{basic_construc}
For an arbitrary  $0 < \delta < 1$,
$0 < \eta < 3/8$ and 
$0 < \ep < 2^{-5} \eta$, set  
$\gamma = 2^{-5}\delta \ep \eta $, and
for $n \in \Nn$ set $N = (1+\ep )n$.
Let $E = (\Rn{N}, \|\cdot\|)$ be an $N$-dimensional Banach space
such that $B_E \subset  B_2^N$.
Let  $\rho \ge 1$,  $0 <{\kappa}  \le 1$ and 
$0 < {a} \le 1$ satisfy
\begin{equation}
 { \vr} (E) \le \rho,  \qquad
{\kappa}  B_1^N \subset B_E,
\qquad  {a} B_2^N \subset B_E.
%%%   \subset  B_2^N.
  \label{spec_pos}
\end{equation}
Then  $E$ admits an $n$-dimensional quotient $F$ such
that  for every operator 
$T: F \to F$ with $T \in \mix_n (\eta n, 1)$,
and every quotient map  $Q: F \to Q(F)$ with 
$\rank\, Q \ge (1 - \gamma) n$,
one has
\begin{equation}
\label{formula_1}
\|QT: F \to Q(F)\| \ge 
c V_{\eta n/4}(E)^{-1} {a}^{ \delta},
\end{equation}
where $c = c(\ep, \eta, {\kappa}, \rho) >0$.
\end{thm}

\noindent\rem
Theorem~\ref{basic_construc} remains valid also for $\de =0$.
In this  case we have $\gamma=0$ 
(note that the constant $c$ {\em does not}
depend on $\gamma$) and the theorem reduces to Theorem 2.2
from \cite{MT3}.

\smallskip
Formula~\ref{formula_1} implies 
that the space $F$ as  well as its further quotients $F_0$
admit  relatively few well bounded operators.
Indeed, one can  formally deduce from it well bounded operators
on these spaces are small perturbations of a multiple
of the identity operator.

\begin{prop}
\label{polmix}
With the same notation as used  in 
Theorem~\ref{basic_construc}
if  $E = (\Rn{N}, \|\cdot\|)$ is an $N$-dimensional 
Banach space such that $B_E \subset  B_2^N$ and 
satisfying (\ref{spec_pos}), then 
$E$ admits an $n$-dimensional quotient $F$ which
satisfies the following two conditions:
\gdef\labelenumi{\theenumi}
\gdef\theenumi{(\roman{enumi})}
\begin{enumerate}
\item  
for every quotient $F_0$ of 
$F$ with $\dim F_0 = k \ge (1 - \gamma)n$  and
every  operator $T: F_0 \to F_0$, 
with $T \in  \mix_k(\eta_0 k ,1)$, 
where $\eta_0 = \eta /(1 - \gamma)$, 
we have
\begin{equation}
\|T\| \ge c V_{\eta n/4}(E)^{-1}  {a}^{ \delta},
\label{quot1}
\end{equation}
\item 
every  Banach space $\tilde{F}$ with
$\dim \tilde{F} =l \le (1+\gamma)n$
such that $F$ is a quotient of $\tilde{F}$
admits a Euclidean norm such that 
every operator $\tilde{T}: \tilde{F} \to \tilde{F}$, 
which is $( \tilde{\eta} l, 1)$-mixing with respect to this norm,
where $\tilde{\eta} = (\eta + 2 \gamma)$, 
satisfies
\begin{equation}
\|\tilde{T}\| \ge c V_{\eta n/4}(E)^{-1}  {a}^{\delta}.
\label{quot2}
\end{equation}
\end{enumerate}
Here $c = c (\ep, \eta,  {\kappa}, \rho)>0$.
\end{prop}
\proof
Let $F$ be the quotient of $E$ satisfying Theorem~\ref{basic_construc}.
In particular,  $F$ admits the natural Euclidean norm inherited from $E$.
Set 
$$
K = \inf_Q \inf \{\|Q\tilde{T}:  F \to Q(F)\| \mid 
\tilde{T}: F \to F, \ \tilde{T} \in \mix_n (\eta n, 1)\},
$$
where the first infimum runs over all 
quotient maps  $Q: F \to Q(F)$ with 
$\rank\, Q \ge (1 - \gamma) n$.
To prove $(i)$, pick a  quotient ${F_0} = F/G_0$, 
with $\dim F_0 = k \geq (1-\gamma)n$.
Identify  $F_0$  with  the  linear subspace $ G_0^\bot \subset F$,
under the norm whose unit ball is  $ Q_{G_0}(B_F)$;
here   $ Q_{G_0}$ is the othogonal projection
with $\ker  Q_{G_0} = G_0$.  
Fix an arbitrary $(\eta' k, 1)$-mixing
operator  $T: F_0 \to F_0$.
Pick any  $\tilde{T}: F \to F$ such that
$Q_{G_0} \tilde{T} = T Q_{G_0}$.
Since  for $x \in G_0^\bot$ we have
$\tilde{T}x = Tx + z$, for some $z \in G_0$, 
then 
$\tilde{T} \in \mix_n(\eta' k, 1)$.
Since $\eta ' k \ge \eta n$, then 
$\tilde{T} \in  \mix_n(\eta n,1)$.
Thus, by  the definition of $K$ we have
%
%\begin{eqnarray*}
$$
\|T: F_0 \to F_0\| = \| Q_{G_0} T Q_{G_0}: F \to F_0\| \\
= \| Q_{G_0} \tilde{T} : F \to F_0\| \ge K,
%   \ge c V_{\eta n/4}(E)^{-1}  {a}^{ \delta}.
%\end{eqnarray*}
$$
and  $(i)$ by the estimate (\ref{formula_1}).

The proof of $(ii)$ is very similar.
Fix an $l$-dimensional space  $\tilde{F}$
with $l \le (1+\gamma)n$
such that $q_F: \tilde{F} \to F$ is
the quotient map. Consider an arbitrary
Euclidean norm on $\tilde{F}$, with the
unit ball $\tilde{B}_2$, such that
$q_F (\tilde{B}_2)$ is the natural
Euclidean ball on $F$. 
Let  $\tilde{T}: \tilde{F} \to \tilde{F}$ 
be   $(\tilde{\eta} l,1)$-mixing.
If $T: F \to F$ satisfies 
$T q_F = q_F \tilde {T}$ then 
clearly,  $\|\tilde{T}\|\ge \|q_F\tilde{T}\| =\|T\|$.
Moreover, since $\dim \ker q_F \le \gamma n$, then
$T \in  \Mix_n(\tilde{\eta}l - 2 \gamma  n,1)
   \subset \Mix_n(\eta n,1)$.
Thus $\|T\|\ge K$ and 
the lower estimate for $\|\tilde{T}\|$
follows.
\qed

It is well known that if for a space $F$ all mixing
operators have large norms, then  $F$ itself  and other related spaces
have  several structural invariants,
such as  basis constant or symmetry constant
or complexification  constant, also bounded below
(\cf \eg\/ \cite{MT3} Section 6).
We give an example of an estimate of this type.

\begin{cor}
  \label{bas_const}
For an arbitrary  $0 < \ep < 2^{-13}$ and  $0 < \delta < 1$,
set $\gamma =  2^{-13}  \ep\delta  $,
and for $n \in \Nn$ set $m = (1+2\ep )n$.
Let $E$ ba an $m$-dimensional Banach space and  
let $\di_E = \di (E, l_2^m)$.
% Let $|\cdot|_2$ be a Euclidean norm on $E$ constructed
% in Lemma~\ref{quotients} $(i)$.
There  exists  an $n$-dimensional quotient $F$ of $E$ such that 
for any  further quotient $F_0$ of $F$ with
$\dim F_0 \ge (1 - \gamma) n$, and for any   space
$\tilde{F}$  with $\dim \tilde{F} \le (1 + \gamma) n$,
for  whose $F$ is  a quotient,  denoting by $F'$ either $F_0$
or $\tilde{F}$, we have
\begin{equation}
 \bc (F')  \ge c V_{2^{-10} n}(E)^{-1}  \di_E^{-\delta}.
\label{11}
\end{equation}
Moreover, for any Banach space  $Z$ such that every 
$n$-dimensional subspace of $Z$ is $D$-Euclidean we have
\begin{equation}
 \bc (F' \oplus_2 Z)  \ge   c D^{-1}  V_{2^{-10} n}(E)^{-1/2}  
                \di_E^{-\delta/2}.
  \label{12}
\end{equation}
Here $c = c(\ep, \delta) >0$.
\end{cor}

The choice of $\eta = 2^{-8}$ made below when applying 
Proposition~\ref{polmix} was  done for convenience of
past references. 
With appropriate  modifications
the same argument would work
for an arbitrary $ 0 < \eta < 3/8$.

\smallskip
\noindent\proof
By  Corollary~\ref{spec_pos_cor} $(i)$, 
there exists  a $(1+\ep)n$-dimensional
quotient  $\tilde{E}$  of $E$ and a Euclidean ball
on $\tilde{E}$, $\tilde{B}_2 \supset B_{\tilde{E}}$,
such that  (\ref{spec_pos})  holds, with
some $\rho = \rho(\ep)$, $ {\kappa}=  {\kappa}(\ep)$
and $ {a} = (2^{1/2}\di_E)^{-1}$.

Fix $\eta = 2^{-8}$.
Then  operators on $F'$ which are
$(5 \cdot 2^{-9} \dim F', 1)$-mixing  are  also
$(\eta' \dim F', 1)$-mixing, 
where $\eta' = \eta_0$ or 
$\eta' = \tilde{\eta}$, depending on
the choice of $F'$ being $F_0$ or $\tilde{F}$.
Therefore, 
by  Proposition~\ref{polmix}, these operators
have norms bounded below
by $K =  c V_{ 2^{-10} n}(E)^{-1}  \di_E^{-\delta}$. 
Thus the conclusion
follows from Theorem 2.1  in \cite{MT3}.
\qed

\noindent\rem
The quotient space $F$ itself satisfies  (\ref{11}) and (\ref{12})
as well, with  $\delta =0$. This  was the  content of Theorem 2.4 in 
\cite{MT3}, and it  followed  from the  construction in \cite{MT3},
Theorem 2.2, which preceded   the present construction.

\medskip
It is of independent interest to study a relationship
between various  $s$-numbers of operators acting in
spaces discussed in Theorem~\ref{basic_construc}.
The advantage of this approach lies in the fact that
resulting estimates  are valid for all operators,
and not only for mixing ones.
%% do not refer to any Euclidean structure
%% on the space.

Let us recall  relevant definitions. Let $X$ and $Y$ be
Banach spaces and let $T: X \to Y$ be a bounded operator.
Let $k$ be a positive integer. 
The $k$th Kolmogorov number $d_k(T)$ is defined by
$$
d_k(T) = \inf_{Z \subset Y} \sup_{x \in B_X} \inf_{y \in Z}
\|T x - y \|,
$$
where the infimum on $Z$ runs over all subspaces $Z$ of $Y$
with $\dim Z < k$.
The dual concept is that of Gelfand
numbers which are defined by
$$
c_k(T) = \inf \{ \|T|_Z \| \mid Z \subset X, \codim Z < k  \}.
$$
We have $c_k(T) = d_k (T^*)$  for arbitrary $X$ and $Y$ and $T$.

As a consequence of Theorem~\ref{basic_construc} we get.

\begin{prop}
\label{s_numbers}
For an arbitrary  $0 < \ep < 2^{-10}$ and $0 < \delta < 1$
set $\gamma = 2^{-10}\delta \ep$, and
for $n \in \Nn$ set $N = 2(1+ 2\ep )n$.
Let $E = (\Rn{N}, \|\cdot\|)$ be an $N$-dimensional Banach space.
Then  $E$ admits an $n$-dimensional quotient $F$ 
such that  for every operator
$T: F \to F$ we have
$$
d_{\gamma n}(T) \ge c V_{ 2^{-7}n }(E)^{-1} \di_E^{-\delta}
      \inf_{\lambda \in \Rn{}} c_{n/4}(T - \lambda \Id_F),
$$
where $c= c(\ep) > 0$.
\end{prop}

\proof
Let $F_0$ be a $(1+\ep)n$-dimensional quotient  of  $E$
satifying condition $(ii)$ of Corollary~\ref{spec_pos_cor},
for   a  Euclidean ball $B_2$ and  some
$ 0 < {\kappa} =  {\kappa} (\ep)  \le 1$,
$\rho= \rho(\ep)$
and $ 0 < {a} \le 1$,
Let $F$ be a quotient of $F_0$ satisfying  (\ref{formula_1}) of
Theorem~\ref{basic_construc},
with $\eta = 2^{-5}$.
Let $Q : F_0 \to F$ be the quotient
map.  Let be the Euclidean norm on $F$
corresponding to $Q(B_2)$; in particular,
$\|x\| \le   {a}^{-1}  \|x\|_2$
for $x \in F$. Then  (\ref{prop_vol}) implies
that 
$$
\bigl(\vol B_F/ \vol Q( {a} B_2^N) \bigr)^{1/n} 
\le C(\ep) \rho.
$$
By the volume ratio argument,  pick  a $15n /16$-dimensional 
subspace  $H$ of $F$  such that 
\begin{equation}
\label{sub_vol}
 {A a}^{-1} \|x\|_2 \le \|x\| \le {a}^{-1}  \|x\|_2
\qquad \mbox{for} \quad x \in H,
\end{equation}
with   $A = A(\ep) = (C(\ep) \rho)^{16} $.

Fix  an arbitrary $T: F \to F$  satisfying
$$
 \inf \{\|(T - \lambda \Id)|_G \|_2
\mid {\lambda \in \Rn{}}, G \subset F, \dim G = 7 n / 8 \}=1.
$$
Pick $G_0 \subset F$ with $\dim G_0 = 7 n / 8 $ 
and $\lambda_0$ such that
$\|(T - \lambda_0 \Id)|_{G_0} \|_2 \le 2$.
Set $G_1 = \{x \in H \cap G_0 \mid Tx \in H\}$.
Then $\dim G_1 \ge 3 n /4$. 
By (\ref{sub_vol}) one has
$$
\|(T - \lambda_0 \Id)x\| \le 2 A \|x\|
\qquad \mbox{for} \quad x \in G_1.
$$
Hence $\|(T - \lambda_0 \Id)|_{G_1}\| \le 2 A $,
which means
\begin{equation}
\label{upper}
  \inf_{\lambda \in \Rn{}} c_{n/4}(T - \lambda \Id_F)
       \le 2 A.
\end{equation}

On the other hand,  the normalization  condition for  $T$ 
and Lemma 2.1 in \cite{M2} imply
that $8\, T$ is $(2^{-5}n, 1)$-mixing.
Denoting the right hand side of (\ref{formula_1}) by $K$,
we get 
$$
\inf \|QT: F \to Q(F)\| \ge K/8,
$$
with the infimum
taken over all  quotient  maps $Q: F \to Q(F)$  with 
$\rank Q \ge (1 - \gamma) n$. 
By the definition of Kolmogorov numbers
this means that $d_{\gamma n} (T) \ge K/8$.
Combining this estimate with (\ref{upper}) we conclude the proof.
\qed
\section{Characterizations of infinite-dimensional\protect\\ spaces}
\label{charact}

As mentioned in the introduction,  if a Banach
space $X$ has the property that its all subspaces have a basis
with a uniform upper bound for the basis constant, then $X$ is
of weak cotype 2. In fact, if there exists  $M < \infty$
such that  $\bc (E) \le M$  for every subspace $E$ of $X$,
then $ wC_2(X) $ admits an upper estimate 
by a function  of $M$.
Indeed, let $F$ be an arbitrary finite-dimensional
quotient  of $X^*$ and consider the Euclidean structure
of $F$ determined by the ellipsoid of minimal volume
containing the unit ball $B_F$. Then 
the  dual version of Theorem 2.4 in \cite{MT3} implies that
$V_{\beta n}(F) \ge c/M$,   for some universal constants
$\beta >0$ and $c >0$.
Then the conclusion follows immediately
from  Lemma~\ref{volumes}, by passing back to the space $X$.

A similar general  line of argument 
is used to prove characterizations of weak cotype 2 spaces
in terms of the basis constant; 
we also obtain related   characterizations 
in terms of an  existence of uniformly bounded projections
and mixing operators.

\begin{thm}
\label{rwaga}
There exists a constant $\gamma_0 >0$ such that a
Banach space $X$ is of weak cotype 2 if and only if
there exist a constant $M \ge 1$ such that 
every finite-dimensional
subspace $E$ of $X$ contains a subspace $E_0 \subset E$ with
$k=\dim E_0 \ge (1 - \gamma_0) \dim E$
satisfying one of the following conditions:
\gdef\labelenumi{\theenumi}
\gdef\theenumi{(\roman{enumi})}
\begin{enumerate}
\item  $\bc (E_0) \le M$,
\item there exists a projection $Q: E_0 \to E_0$ of
rank $k/8$ such that $\|Q\| \le M$,
\item for every Euclidean norm on $E_0$ there exists 
an operator $T: E_0 \to E_0$ which is
$(k/8, 1)$-mixing with respect to this norm,
such that $\|T\| \le M$.
\end{enumerate}
Moreover, if one of the conditions (i)--(iii) holds then 
$ wC_2(X) \le C M^4$,
where $C$ is a universal constant.
\end{thm}

The counterpart of this result for superspaces is  less satisfactory,
as it gives  implication in  one direction only.

\begin{thm}
\label{nwaga}
There exists a constant $\gamma_0 >0$ such that whenever
$X$ is a Banach space for which
there exist a constant $M \ge 1$ such that for
every finite-dimensional subspace $E$ of $X$ there 
is a Banach space $\tilde{E} \supset E$ with
$k = \dim \tilde{E} \le (1 + \gamma_0) \dim E$
satisfying one of the following 
conditions
\gdef\labelenumi{\theenumi}
\gdef\theenumi{(\roman{enumi})}
\begin{enumerate}
\item  $\bc (\tilde{E}) \le M$,
\item  there exists a projection $Q: E_0 \to E_0$ of
rank $k/8$ such that $\|Q\| \le M$,
\item for every  Euclidean norm on $\tilde{E} $ there exists 
an operator $T: \tilde{E} \to \tilde{E}$ which is
$(k/8, 1)$-mixing with respect to this norm,
such that $\|T\| \le M$,
\end{enumerate}
then $X$ is of  weak  cotype 2 and 
$ wC_2(X) \le C M^4$,
where $C$ is a universal constant.
\end{thm}

\noindent\rem
Beside conditions $(i)$ and $(ii)$ of the above theorem there is a
number of other invariants whose uniform boundedness implies condition
$(iii)$. They are  symmetry constant (\cf\/ \cite{M}, \cite{M3}),
complexification constant (for real Banach spaces) or the
Banach--Mazur distance from the space to its complex conjugate (for
complex Banach spaces) (\cf\/ \cite{S2}). All these 
and other invariants could be used for 
versions of  all theorems of this section
(\cf\/ also \cite{MT3}, Section 6).

\smallskip
\noindent\rem
Note that there is an essential difference betwen the Euclidean case
in  the definition of  weak cotype 2
and the results  above. Namely,  we do
not know whether analogous  characterizations  are valid
for an {\em arbitrary} proportion $\de  \in (0,1)$
(not necessary $\de > 1 - \gamma_0$). Recall that
it is so in  the Euclidean case: 
if for a Banach space $X$  there exists  $\de_0 >0$ such that 
every finite dimensional subspace $E$ of $X$ contains 
a $C_0$-Euclidean subspace 
$E_0 \subset E$ with $\dim E_0 \ge \de_0 \dim E$,
then an analogous condition holds for every $\de \in (0,1)$,
with the constant $C$ depending on $\de$.
It seems that the present difficulty
is connected with a problem of Pe{\l}czy\'nski in \cite{Pe}
whether every finite dimensional Banach space $E$ can be embeded into
a Banach space $F$ with $ \dim F \le 2 \dim E$ having 
a nice Schauder basis. For $\dim F$ close enough to $\dim E$,
this question was answered in \cite{MT2} in the negative.

\bigskip

\noindent {\bf Proof of Theorem~\ref{rwaga}} 
Clearly, the weak cotype 2 assumption implies property $(i)$, which
implies $(ii)$, which implies $(iii)$.
We shall prove that conversely, property  $(iii)$ implies 
that $X$ is of weak cotype 2. 
Set $\gamma_0 = 2^{-20}$. Under our assumptions we have 
the following.

\medskip
\noindent {\bf Claim}
{\it For every finite-dimensional subspace $E \subset X$ there exists
a subspace $H \subset E$ with $\dim H \ge 2^{-9} \dim E$ and 
$\di_H \le c M^2 \di_E^{1/2}$, where $c$ is a universal constant.}

\medskip
\noindent {\bf Proof of the Claim}
Let $\ep = 2^{-9}$, $\eta = 2^{-4}$, $ \delta = 2^{-2}$, and 
$\gamma = \gamma_0$. Fix an arbitrary 
finite-dimensional $E \subset X$ and let $n \in \Nn$
be such that $\dim E = (1 + 2\ep) n$.

Set $Z = E^*$.
Let $Z_1$ be a quotient of $Z$ with 
$\dim Z_1 = (1+\ep)n$ 
satisfying  condition  $(i)$ of Corollary~\ref{spec_pos_cor}.
We also  fix
the Euclidean structure on $Z_1$ introduced in this condition.
Let $Z_2$,   with  $\dim Z_2 = n$,  be a quotient of $Z_1$
constructed in Proposition~\ref{polmix}.
Now we use  the dual form of $(iii)$ 
valid for every finite-dimensional quotient of $X^*$
(\cf\ Proposition~\ref{mix} $(iv)$).
It follows that  $Z_2$  has a
quotient $Z_3$ with $\dim Z_3 = k \ge (1-\gamma) n$
which admits  a $(k/8, 1)$-mixing operator  $S_0$
with $\|S_0\| \le M$. 
On the other hand, 
by the choice of $Z_2$,  every
$(\eta'k,1)$-mixing operator $T$ on $Z_3$ satisfies 
$$
\|T\| \geq cV^{-1}_{2^{-6} n}(Z_1)\di^{-\de}_{Z},
$$
where $\eta'=\eta /(1-\gamma)$. Since
$\eta' \le 1/8$, the same estimate holds for $S_0$.
Thus
$V_{2^{-6} n}(Z_1) \geq c \di^{-\de}_{Z}M^{-1}$.
Applying Lemma~\ref{volumes} with
$\beta = 2^{-6}$ and $ \sigma = 1$,
we obtain a quotient $Z_4$ of $Z_3$ with 
$\dim Z_4 \geq 2^{-7}k \geq 2^{-9} \dim E $ and 
$\di_{Z_4} \leq C\di^{2\de}_{Z}M^2$.
Since $\di_{Z} = \di_E$, the proof of the Claim is 
concluded by setting $H=Z^*_4$.

\smallskip
Passing to the proof of the theorem,  
fix an arbitrary finite-dimensional subspace
$X_0 \subset X$.
Denote by $D$  the smallest number
such that every subspace $X_1$ of $X_0$ contains a
subspace $G \subset X_1$ with $\dim G \ge  \dim X_1/2$ and
$\di_G \le D$. Clearly, $D$ is finite.
% by (\ref{weak_cot_2}),
% it is enough to show that  $(iii)$ implies that
% $D \le C_1 M^4$ for some numerical constant $C_1$. 

We then know  that for any $m$ and any subspace $X_1$ of $X_0$ with
$\dim X_1 = m$ there is a subspace $H \subset X_1$ with
$\dim H \ge 2^{-10} m$ such that $\di_H \le c M^2 D^{1/2}$.
Indeed, first  pick $E \subset X_1$ with
$\dim E \ge m/2$  such that $\di_E \le D$, and
then apply Claim to obtain $H$.
Specifying  $m = (1/2) \dim X_0$ we can 
use Lemma~\ref{bighilbert} with  $\xi = 1/2$,
$\delta = (1 - 2^{-10}) \xi$ 
and $\eta = 2^{-10}\xi$, to  get a subspace
$ \tilde{H} \subset X_0$ satisfying 
$\dim  \tilde{H} \ge (1/2)\dim X_0$
and $\di_{\tilde{H}}\le c' M^2 D^{1/2}$,
where $c'$ is a universal constant.
By the definition of $D$, this implies
$D \le c' M^2 D^{1/2}$, hence 
$D \le c'' M^4$. By (\ref{weak_cot_2}), this completes
the proof.
\qed

Theorem~\ref{nwaga} has  almost identical proof,
with  condition $(ii)$   of  Proposition~\ref{polmix} 
replacing $(i)$. We shall omit further details.

\noindent \rem
Using Lemma~\ref{volumes} in a more delicate way and choosing $\de$
sufficiently small one can get in the theorems above $wC_2(X) \le
C(\sigma) M^{1+ \sigma}$, for every $\sigma>0$.

\medskip
We now pass to characterizations of weak Hilbert spaces.

\begin{thm}
\label{rweakhil}
There exists a constant $\gamma_0 >0$ such that a
Banach space $X$ is a weak Hilbert space if and only if
there exists a constant $M \ge 1$ such that 
one of the following conditions is satisfied for
every subspace $Y \subset X$:
\gdef\labelenumi{\theenumi}
\gdef\theenumi{(\roman{enumi})}
\begin{enumerate}
\item   every finite-dimensional
subspace  $E$  of $Y$ contains a subspace $E_0$
with $\dim E_0 \ge (1 - \gamma_0) \dim E$
which admits  a projection $Q: E_0 \to E_0$ with
$\rank Q = \dim E_0/8$ and  $\|Q\| \leq M$,
and  \\ every finite-dimensional quotient $F$ of $Y$
admits   a quotient $F_0$ with
$ \dim F_0 \ge (1 - \gamma_0) \dim F$
which admits a projection $R: F_0 \to F_0$
with $\rank R = \dim F_0/8$ and $\|Q\| \leq M$,
\item  for every finite-dimensional
subspace  $E$  of $Y$ there is a Banach space $\tilde{E}$ containing
$E$ with $\dim \tilde{E} \le (1 + \gamma_0) \dim E$
which admits a  projection $\tilde{Q}: \tilde{E} \to \tilde{E}$
with $\rank \tilde{Q} =\dim \tilde{E}/8$  
and $\|\tilde{Q}\| \le M$,
and \\ for every finite-dimensional quotient $F$ of $Y$
there is a Banach space $\tilde{F}$ which has $F$ as a a quotient
and $ \dim \tilde{F} \le (1 + \gamma_0) \dim F$
which admits a  projection $R: \tilde{F} \to \tilde{F}$
with $\rank R = \dim \tilde{F}/8$ and $\|R\| \leq M$.
\end{enumerate}
\end{thm}
\proof
Clearly, if $X$ is a weak Hilbert space then both $(i)$ and $(ii)$
are satisfied. Conversely, assuming that one of
conditions $(i)$ or $(ii)$ is satisfied,
by Theorem 1 in \cite{MT2}, $X$ does not contain $l_1^n$'s
uniformly. Hence, by Pisier's result,  \cite{P1},
(\cfalso\ \cite{P}, Theorem 11.3), $X$ is $K$-convex.

If $(i)$ holds, by Theorem~\ref{rwaga} we infer that
both $X$ and $X^*$ are of weak cotype 2. Thus $X$ is a weak
Hilbert space. In case of $(ii)$, we use Theorem~\ref{nwaga}.
\qed

\noindent\rem
If $X$ has the approximation property, the above theorem  remains
valid if we restrict ourselves to the case $Y = X$ only.

\noindent\rem
In Theorem~\ref{rweakhil}, the assumption on the existence of uniformly
bounded rank $k/8$ projections can be replaced by the existence of
uniformly bounded $(k/8,1)$-mixing operators or by a uniform bound for
basis constants.
\section{Subspaces and quotients of proportional dimension}
\label{fin_charac}

The main result of this section is concerned with
subspaces and quotients of a fixed  proportional dimension.
\begin{thm}
\label{main2}
Let $0 < \al <(1+ 2^{-8})^{-1}$ and let $M \ge 1$.
Let $G$ be an $n$-dimensional Banach space.
If every $\al n$-dimensional subspace  $E$
and every $\al n$-dimensional quotient $F$ of
$G$ have the  basis constants 
$\bc (E) \le M$ and $\bc (F) \le M$,
then $G$ is a weak Hilbert space  and the
weak type 2 and the  weak cotype 2 constants satisfy
the estimates
$$
wT_2(G) \le C(\al)M^{140 /3}  \quad \mbox{and}\quad  wC_2(G) \le 
 C(\al)M^{100 /3}.
$$
\end{thm}

The proof of the theorem requires several steps.
To begin with we consider only one-sided assumptions
on $G$, that is, the assumptions on its subspaces.
In such a situation,  the following proposition 
establishes a weak cotype 2 property,
provided that the space admits a nice
direct sum decomposition.

\begin{prop}
\label{cotype}
Let $0 < \al < 1$,  $M \ge 1$
and let  $G$ be an $n$-dimensional Banach space.
Assume that $G = Z \oplus_2 G_0$, with $\dim Z \ge \al n$.
If every $\al n$-dimensional subspace $E$  of $G$ has the  basis
constant $\bc (E) \le M$, 
then the  weak cotype 2 constant of $G_0$
satisfies  $wC_2(G_0) \le C(\al) (M \inf_H \di_H )^4$,
where the infimum runs over  all
$\al n$-dimensional subspaces  $H \subset Z$.
\end{prop}
\proof
Let  $H \subset Z$ be an $\al n$-dimensional subspace.
We shall show that every subspace $G_1$ of $G_0$
contains a subspace $\tilde{E}$, with 
$\dim \tilde{E} \ge    2^{-12} \min(\dim G_1, \al n)$
such that  $\di_{\tilde{E}} \le    C (M \di_H )^4$.
By (\ref{weak_cot_2}),
this  will imply 
$wC_2(G_0) \le C(\al) (M  \di_H )^4$,
hence the conclusion will follow by passing 
to the infimum over  $H$.
Obviously it is enough to consider
the case $\dim G_1 \le  \al n$ only.

Fix   an arbitrary  subspace $G_1$  of
$G_0$ with $k = \dim G_1 \le  \al n$
and let $F = G_1^*$.  
By the remark following  Corollary~\ref{bas_const} 
there exists a quotient $F_1 $ of $F$ with
$ k' = \dim F_1 \ge (1+2^{-12})^{-1}k $ 
satisfying (\ref{12}) with  $\delta =0$.
Fix  an arbitrary $H_1 \subset H$ with 
$\dim H_1 = \al n - k'$. We have
$\di_{H_1^*}= \di_{H_1} \le \di_H$
and  $H_1^*$ is a quotient of $Z$.
Thus
$$
 \bc (F_1 \oplus_2 H_1^*)  
    \ge   c \di_H^{-1}  V_{2^{-10} k'}(F)^{-1/2}.
$$
Observe that $(F_1 \oplus_2 H_1^*)^*$ is an $\al n$-dimensional
subspace of $G$. Hence $\bc (F_1 \oplus_2 H_1^*) = 
\bc ((F_1 \oplus_2 H_1^*)^*)\le M$. Combining the last
two estimates we get
$$
 V_{2^{-10} k'}(F) \ge c' (M \di_H)^{-2}.
$$
By Lemma~\ref{volumes} (with $\sigma = 1$) 
we obtain a  Euclidean quotient
of $F$ of dimension $2^{-11}k' \ge 2^{-12} k$ 
and so  we complete the proof by passing to the dual.
\qed

Let us note that the one-sided boundness assumption alone
still yields the existence of some Euclidean subspaces, but
this time on a proportional level only.
\begin{lemma}
\label{homog}
Let $0 < \al  <(1+ 2^{-12})^{-1}$ and let $M \ge 1$.
Let $G$ be an $n$-dimensional Banach space such that
every $\al n$-dimensional subspace  $E$ of $G$ has 
the  basis constant  $\bc (E) \le M$.
Then  for every $ (1 + 2^{-12})\al < \la \le 1$,
every $\la n$-dimensional subspace $G_0$ of $G$
contains a subspace $H$ with
$\dim H = (\la - \al) n \ge 2^{-13}\la n$ such that 
$\di_H \le C(\al, \la)  M^{2}$. 
\end{lemma}
\proof
A similar argument as in the proposition above,
in which  the use of (\ref{12})  is replaced
by (\ref{11}), shows that  
every subspace $G_1$ of $G_0$ with $\dim G_1 = \xi n
= (1 + 2^{-12})\al n$ 
contains a $C M^2$-Euclidean subspace
of dimension $2^{-11} \al n$,
where $C$ is a universal constant.
The proof is then concluded by  applying
Lemma~\ref{bighilbert}, with 
$\delta  =  (1-2^{-12})\al $ and $\eta =  2^{-12}\al $,
to any $\la n$-dimensional subspace $G_0$ of $G$.
\qed

An {\sl a priori} argument which we are going to use
is based on a finite-dimensional
version of one of properties characterizing
weak type 2 spaces.
A known argument (\cf\ \cite{P}, Chapter 11),
localized to a fixed proportional-dimensional
level, gives  the following lemma.
Recall that for a Banach space $G$,  $K(G)$ stands for 
the $K$-convexity constant of $G$.
\begin{lemma}
\label{lemma_proj}
Let $0 <\de < \be <1$, $D \ge 1$, and 
let $G$ be an $n$-dimensional Banach space.
Assume that every $\be n$-dimensional
subspace $F$ of $G^*$
contains a subspace $F_1$ with $\dim {F_1} = \de n$
and $\di_{F_1} \le D$. Then for every  $\be n$-dimensional 
subspace $E$ of $G$
there exist
a subspace $H$ with $\dim H = \de n/2$
and a projection $Q: G \to H$ such that 
$\|Q\| \le C(\be, \de) K(G) D \di_E $.
%where $C$ is a universal constant.
\end{lemma}
\proof
Let  $E$ be a  $\be n$-dimensional  subspace of $G$
and let $w: E \to l_2^{\be n}$ be 
an isomorphism  such that
$\|w\| \, \|w^{-1}\|= \di_E$. 
We will show that
there exist an orthogonal rank $(\de/2)n$ 
projection $P$ in $l_2^{\be n}$
and an operator $\tilde{w}: G \to  l_2^{\be n}$,
such that $\tilde{w}= P w$ and $\|\tilde{w}\|
\le C(\be, \de) D K(G) \|w\|$.
Then $H = w^{-1} P (l_2^{\be n})$ and $Q = 
w^{-1} P \tilde{w}$ will satisfy the requirements
of the lemma.

The argument
requires the definition of the $l$-norm of an operator 
$u: l_2^k \to Y$, for any Banach space $Y$, which is provided
e.g. in \cite{P}, Chapter 3. 
Similarly as in the proof of Theorem 11.6 in \cite{P}, 
consider the operator $w^*:  l_2^{\be n} \to G^*/E^{\bot} $.
By Lemma 11.7 in \cite{P}, there exists 
$\tilde{v}:  l_2^{\be n} \to G^*$ such 
that $q  \tilde{v} = w^*$, where $q: G^*  \to G^*/E^{\bot} $
it the quotient map, and that 
$l(\tilde{v}) \le 2 K(G) l(w^*) \le 2 (\be n)^{1/2}K(G) \|w^*\|$.
Note that $ \tilde{v}$ is one-to-one.
Choose 
a subspace $F_1$ of   $ \tilde{v}( l_2^{\be n})$
with $\dim {F_1} = \de n$
and $\di_{F_1} \le D$ and let $F_2 =  \tilde{v}^{-1}(F_1)$.
By well-known properties of operators
acting in Hilbert spaces (\cf\eg\/ \cite{P} Proposition 3.13),
there is a $(\de/2) n$-dimensional subspace $F_3 \subset F_2$
and an orthogonal projection $P$ onto $F_3$
such that 
$$
\|\tilde{v} P  \|\le (\de n /2)^{-1/2} D\, l(\tilde{v})
\le 2^{3/2} (\be / \de)^{-1/2} K(G) D \|w\|.
$$
Therefore the required operator is $\tilde{w}= P \tilde{v}^*$.
\qed

We are finally ready for the proof of the theorem.
\medskip

\noindent{\bf Proof of Theorem~\ref{main2}\ \ }
Fix $0 < \al < 1/4$ and set $\al' = (1 + 2^{-12})\al$.
By Lemma~\ref{homog} with $\la =1$, pick a subspace $H$ of $G$
with $\dim H \ge (1 - \al)n$  and $\di_H \le C(\al) M^2$.
Since all $\al n$-dimensional quotients $F$ of $G$ satisfy
$\bc (F) \le M$, applying the same lemma for $G^*$
and  $\la= (1 - \al) > \al'$, it follows 
that $G^*$ satisfies the assumptions of  Lemma~\ref{lemma_proj}
for $\be = 1 - \al $, $\de = 1 - 2\al $ and $D = C(\al) M^2$.
Therefore there exist  a subspace $H_0$ of $H$
and a projection $Q$ from $G$ onto $H_0$ such that
$\dim H_0 = (1/2 - \al )n$ and $\|Q\| \le C'(\al) K(G) M^4$. 
Notice that $\dim H_0 \ge \al n$.

Set $G_0 = \ker Q$ and apply Proposition~\ref{cotype} to both
$G$ and $G^*$ to conclude that the their weak cotype 2 constants  
satisfy
\begin{equation}
\max\,(wC_2(G), wC_2(G^*)) \le C''(\al) \|Q\| (\di_{H_0} M^2)^4 
      \le C''(\al) K(G) M^{20}.
\label{const}
\end{equation}
In particular, by (\ref{weak_typ_2}), we also have
\begin{equation}
wT_{2}(G) \le K(G) wC_2(G^*) \le
 C''(\al) K(G)^2 M^{20}.
\label{wtype}
\end{equation}

Now we use the result of Pisier stated in 
(\ref{kconvex}), for \eg\  $\theta = 1/6$, to  get
$ K(G) \le  C'''(\al) K(G)^{1/2} M^{40/6}$. 
Thus 
$$
K(G) \le C_0(\al) M^{40/3}.
$$
The proof is then concluded by 
combining this inequality with (\ref{const}) and 
(\ref{wtype}).
\qed

\noindent\rem
Applying Lemma~\ref{homog} with $\sigma$ arbitrarily close to 1
and  Corollary~\ref{bas_const}  with $\ep$ aufficiently small,
one can  replace  the constant $1/4$ in Theorem~\ref{main2}
by any number smaller than 1.
Similar  rearrangements yield
more civilized powers of $M$.
\section{Random quotients}
\label{random}

Fix a  probability space $(\Omega,{\bf P})$ and
let $g_1, \ldots, g_{\ep n}$ be independent standard Gaussian
vectors in $\Rn{n}$ with the
density
$  ({n}/{2\pi})^{n/2} {\rm e}^{-n\|x\|_2^2 /2} $,
with respect to the standard Lebesgue measure in $\Rn{n}$.

For $\om \in \Om$,
define a Gaussian projection
$Q_{\om}: \Rn{N} \rightarrow \Rn{n}$ by
$$
Q_{\omega} (e_i) = \left\{  \begin{array}{ll}
e_i & \mbox{ for $i=1,2,\ldots,n$}\\[2mm]
g_{i-n}(\om) & \mbox{ for $i=n+1,n+2,\ldots,N$.}
\end{array}
\right.
$$

In the theorem below, we denote by
$ G_{\gamma n, n}$  the set of all $\gamma n$-dimensional
subspaces of $\Rn{n}$. For $ H \in  G_{\gamma n, n}$ we denote
by $Q_H$ the orthogonal projection with $\ker Q_H = H$.
\begin{thm}
\label{zmiara}
For an arbitrary  $0 < \delta < 1$,
$0 < \eta < 3/8$ and 
$0 < \ep < 2^{-5} \eta$, set  
$\gamma = 2^{-5}\delta \ep \eta $, and
for $n \in \Nn$ set $N = (1+\ep )n$.
Let $E = (\Rn{N}, \|\cdot\|)$ be an $N$-dimensional Banach space
such that $B_E \subset  B_2^N$.
Let  $\rho \ge 1$  and 
$0 < {a} \le 1$ satisfy
$$
 { \vr} (E) \le \rho
\qquad  {a} B_2^N \subset B_E.
%%%   \subset  B_2^N.
$$
% 
% For an arbitrary $0 < \ep' < 2^{-5}$, $0 < \delta < 1$
% and $0 < \eta < 3/8$,
% set $\ep = \ep' \eta$ and  
% $\gamma = 2^{-5}\delta \ep \eta $, and
% for $n \in \Nn$ set $N = (1+\ep )n$.
% Let $E = (\Rn{N}, \|\cdot\|)$ be an $N$-dimensional Banach space.
% Let  $\rho \ge 1$  and  $\overline{a} \le 1$ satisfy
% $$
%  { \vr} (E)= \rho,  \qquad
% \qquad  \overline{a} B_2^N \subset B_E \subset  B_2^N.
% $$
%
There exists $0 < c = c(\ep,\eta, \rho) < 1$ such that if
$\tilde{\Om}$ denotes the set
\begin{eqnarray*}
\tilde{\Om} = \{\om \in \Om &\mid& 
\|Q_H T: Q_{\om}(l_1^N) \to Q_H Q_{\om}(E)\|
\ge c   V_{\eta n/4 }(E)^{-1}  {a}^{\delta},\\
&& \mbox{for every } T \in {\Mix_n}(\eta n,1)
\mbox{ and every }
H \in G_{\gamma n, n} \},
\end{eqnarray*}
then ${\bf P} (\tilde{\Om}) \ge  1 - c_1^{n^2}$,
where  $0 < c_1 <1$ is an absolute constant.
\end{thm}

To deduce Theorem~\ref{basic_construc} pick  
$\om \in \tilde{\Om}$ and set $F = Q_{\om}(E)$.
Now it is enough to observe
that since $ {\kappa}  B_1^N \subset B_E$  then
$$
\|Q_H T: Q_{\om}(E) \to Q_H Q_{\om}(E)\|\ge  {\kappa}
\|Q_H T: Q_{\om}(l_1^N) \to Q_H Q_{\om}(E)\|.
$$

For   $E = l_1^N$, Theorem~\ref{zmiara} was proved
recently in \cite{MT2} Theorem 4. In the general case
the argument follows the steps from  \cite{MT2} blended with a
technique which enables to pass from quotients of $l_1^N$
to quotients of arbitrary Banach spaces, as presented
in \cite{MT3}, Section 5. Therefore we shall only briefly
discuss the main points, referring the reader to
\cite{MT2} and \cite{MT3} for the details.

Passing to the description of the proof of
Theorem~\ref{zmiara}, we require additional notation.
For every $\om \in \Om$ let
$H_{\om}= {\rm span}\,[g_1(\om), g_2(\om),\ldots,g_{\ep n}(\om)]$.
If $H \in G_{\gamma n,n}$, 
let $Q_{\om,H}$ be the orthogonal projection in $\Rn{n}$ with
$\ker Q_{\om,H} = H + H_{\om}$.

Let
$$
\Om_0 = \{ \om \in \Om \mid 1/2 \le \|g_i(\om)\|_2 \le 2
\quad \mbox{for all } i = 1, \ldots, \ep n \}.
$$
Fix $T \in {\Mix_n}(2\eta n/3,1)$. By the definition
of the mixing class, there is $G \subset \Rn{n}$, $\dim G = 2\eta n/3$
such that $\|P_{G^{\perp}}T x\|_2 \geq \|x\|_2$ for every $x \in G$.
The well-known  argument
on half-dimensional circular sections of
an ellipsoid yields  that there exists $G_0 \subset G$
with $\dim G_0 = \eta n/3$
and $\lambda \ge  1$ such that 
$\|P_{G^{\perp}}T x\|_2 = \lambda \|x\|_2$
for every $x \in G_0$. 
For every $\om \in \Om$ and $H \in G_{\gamma n,n}$ fix
an orthogonal projection $Q_{\om, H, G_0}$ in $\Rn{n}$ with
$$
\ker Q_{\om, H, G_0} \supset H_{\om} + H + G + 
P_{G^{\perp}} T P_{G_0 ^{\perp}}(H_{\om}),
$$
and $\rank Q_{\om, H, G_0} = \eta n /4$.

Set
\begin{eqnarray}
\Omega_{T, H} =
\{\om \in \Om_0 &\mid& Q_{\om, H, G_0} T P_{G_0} g_j  \in 4 \lambda 
    \alpha \sqrt{\eta}  {a}^{\delta} V_{\eta n/4}^{-1}
          Q_{\om, H, G_0}Q_\om (B_E)\nonumber\\
&& \mbox{for } j=1,2,\ldots,\ep n\}.
\label{omega_TF}
\end{eqnarray}

\begin{lemma} 
\label{singleTF}
Let $H \in G_{\gamma n,n}$. Then
$$
{\bf P}(\Om_{T, H}) \leq \bigl(C_0\alpha \sqrt{\eta} 
         {a}^{\delta}\bigr)^{\ep \eta n^2/4},
$$
where $C_0 \ge 1$ is an absolute constant.
\end{lemma}
\proof
Set $\widetilde{Q}_{\om}= Q_{\om, H, G_0}$ for $\om \in \Omega$.
For $j=1,2,\ldots,\ep n$ define $g'_j=P_{G_0}g_j$ and
$g''_j=P_{G_0^{\perp}}g_j$.  
Similarly as in \cite{MT2}, Lemma 7,
$\widetilde{Q}_{\om}$ is independent of the $g'_j$'s.
For every fixed $j=1,2,\ldots,\ep n$ we have
\begin{eqnarray*}
\lefteqn{\{\om \in \Omega \mid \widetilde{Q}_{\om} T P_{G_0}g_j 
\in 4 \lambda \alpha 
      \sqrt{\eta}  {a}^{\delta} V_{\eta n/4}^{-1}
\widetilde{Q}_{\om}Q_\om (B_E) \}} \\   
&=& \{\om \in \Omega \mid \widetilde{Q}_{\om} T g'_j 
\in 4 \lambda \alpha  
     \sqrt{\eta}  {a}^{\delta} V_{\eta n/4}^{-1}
              \widetilde{Q}_{\om}Q_\om(B_E)\}.
\end{eqnarray*}

Since $G \subset \ker Q_\om$, the definition of $\lambda$
implies that  $\lambda ^{-1}\widetilde{Q}_{\om}T $
is a contraction in the Euclidean norm on $\Rn{n}$.
Moreover, it
has $k$ $s$-numbers equal to $1$, with
$k \ge \eta n /3 - \gamma n - 2 \ep n \ge \eta n/4$.
Hence, using Claim 6.2 in [Sz.1] (with $n/3$ replaced
by $\eta n/4 $), and the fact that 
$\sqrt{3/ \eta} g'_j$ is a standard Gaussian variable in $G_0$,
(note that $3/ \eta = n/\dim G_0 $ ),
\cfeg [Sz.1] (3.3), we have, by the definition of 
$V_{\eta n/4}$,
\begin{eqnarray*}
\lefteqn{ {\bf P}(\{\om \in \Omega \mid  \widetilde{Q}_{\om} T g'_j 
\in 4 \lambda \alpha 
       \sqrt{\eta}  {a}^{\delta} V_{\eta n/4}^{-1}
               \widetilde{Q}_{\om}Q_\om (B_E )\})} \\
&\leq&
{\bf P}(\{\om \in \Omega  \mid (\lambda^{-1}\widetilde{Q}_{\om} T) 
      (\sqrt{3/\eta}g'_j) 
\in 4 \sqrt{3} \alpha  {a}^{\delta} V_{\eta n/4}^{-1}
               \widetilde{Q}_{\om}Q_\om (B_E)\}) \\
&\leq&
\left(c' 4 \sqrt{3} \alpha  {a}^{\delta} \right)^{\eta n/4}.
\end{eqnarray*}  
where $c'$ is an absolute constant.
Hence
$$
{{\bf P}(\Omega_{T, H})} \leq
(c'\alpha  {a}^{\delta})^{\ep \eta n^2 /4},
$$
which concludes the proof of the lemma.
\qed

The next lemma is a restatement of Lemma 7.3 in [Sz.1].

\begin{lemma}
\label{dual_sud}
For every $0 <\sigma <1 $, 
the set 
$$
{\cal P}_{k,n} = \{ P: \Rn{n} \to \Rn{n} \mid
P \mbox{ an orthogonal projection with } \rank P = n-k\}
$$
admits a $\sigma$-net $\cal M$ in the operator norm
in $l_2^n$ with the cardinality 
$|{\cal M}| \le C^{n^2} \sigma ^{-n k}$,
where $C >1$ is an absolute constant.
\end{lemma}

Using Lemma~\ref{singleTF} and Lemma~\ref{dual_sud} with
$\sigma = \alpha  {a}^{1 + \delta}/4$ and $k = \gamma n$,
the same argument as in the proof of Proposition 5 in \cite{MT2}
yields.

\begin{prop} 
\label{prop_singleTp}
Let $0<\alpha, \delta <1$ and $0 < \eta < 3/8$.
For an operator $T \in {\Mix_n}(2 \eta n/3,1)$
set
\begin{eqnarray} 
\Om_T =  \{\om \in \Om_0 &\mid& \|Q_{\om, H} T : Q_\om (l_1^N) \to
Q_{\om, H}Q_\om (E)\| \leq 2\alpha 
       \sqrt{\eta}  {a}^{\delta} V_{\eta n/4}^{-1}
  \nonumber \\
&& 
\mbox{for some } H \in G_{\gamma n,n}\}.
\label{singleTe'}
\end{eqnarray}  
Then for every $T \in {\Mix_n}(2 \eta n/3,1)$ one has
\begin{equation}   
{\bf P}( \Om_T ) \leq C^{n^2}
       (4/ {a}^2 \alpha)^{\gamma n^2} 
       \bigl(C_0\alpha \sqrt{\eta} 
             {a}^{\delta}\bigr)^{{\ep \eta n^2}/4},            
\label{est_single}
\end{equation}
where $C\ge 1$ and $C_0\le 1$ are absolute constants.
\end{prop}

The rest of the proof of Theorem~\ref{zmiara}
is essentially the same as of Theorem 4 in \cite{MT2}. 
One has to replace
Proposition 6 there by Lemma 5.3 from \cite{MT3},
with $A = 2 \alpha \sqrt{\eta} 
{a}^{\delta} V_{\eta n/4}^{-1}$
and choosing $\alpha > 0$ sufficiently small
as to ensure that
$(4 c(\rho))^{1+{\ep^2}} C\, 
C_0^{\ep \eta /4}(\alpha)^{\ep\eta /8} \le 1/2$.

\address

\begin{thebibliography}{F-L-Mx}
\bibitem[B]{B} {\bf Bourgain, J.}, On finite-dimensional homogeneous
Banach spaces. {\sl GAFA Israel Seminar} 1986-87, 
{\sl Lecture Notes in Math.,} Vol. 1317, Springer, 232--239.
%
\bibitem[B-S]{BS} {\bf Bourgain, J. \& Szarek, S.~J.}, The Banach-Mazur
distance to the cube and the Dvoretzky-Rogers factorization. {\sl
Israel J. of Math.}, 62 (1988), 169-180.
%
% \bibitem[G]{G} {\bf Gluskin, E. D.}, Finite-dimensional analogues of
% spaces without a  basis. {\sl Dokl. Acad. Nauk SSSR}, 261 (1981),
% 1046--1050, (in Russian).
%
% \bibitem[J-P]{JP} {\bf Johnson, W.~B. \& Pisier, G.}, Proportional UAP
% characterizes weak Hilbert spaces. {\sl J. of London Math. Soc.}, 
%
\bibitem[Ma.1]{M} {\bf Mankiewicz, P.}, Finite dimensional Banach spaces
with symmetry constant of order $\sqrt{n}$.
{\sl Studia Math.}, 79, (1988), 193--200.
%
\bibitem[Ma.2]{M3} {\bf Mankiewicz, P.}, Subspace mixing properties
of operators in $\Rn{n}$ with applications to Gluskin spaces.
{\sl Studia Math.}, 88 (1988), 51-67.
%
\bibitem[Ma.3]{M2} {\bf Mankiewicz, P.},  Factoring the identity 
operator on a subspace of $l^n_{\infty}$. {\sl Studia Math.}, 95 (1989),
133-139.
%
\bibitem[M-T.1]{MT1} {\bf Mankiewicz, P.} \&  {\bf Tomczak-Jaegermann, N.},
A solution of the finite-di\-men\-sio\-nal homogeneous Banach spaces
problem. {\sl Israel J. Math.} 75 (1991), 129--159.
%
\bibitem[M-T.2]{MT2} {\bf Mankiewicz, P.} \&  {\bf Tomczak-Jaegermann, N.},
Embedding subspaces of $l^n_{\infty}$ into spaces with Schauder basis.
{\sl Proc. of the AMS} 117 (1993), 459--465.
%
\bibitem[M-T.3]{MT3} {\bf Mankiewicz, P.} \&  {\bf Tomczak-Jaegermann, N.},
Schauder bases in subspaces of quotients of $l_2(X)$. To appear in
{\sl Amer. J. of Math.}
%
\bibitem[Mi]{Mi} {\bf Milman, V.~D.}, In\'egalit\'e de Brunn-Minkowski inverse
et applications \a` le th\'eorie local des espaces norm\'es. {\sl C. R. Acad.
Sci. Paris}, 302 S\'er. 1 (986), 25-28.
%
\bibitem[M-P]{MP} {\bf Milman, V. D.} \& {\bf Pisier, G.},
Banach spaces with weak cotype 2 property. {\sl Israel J. Math.} 54 (1986),
139--158.
%
\bibitem[Pe.1]{Pe} {\bf Pe{\l}czy\'nski, A.}, {\sl Notes in Banach spaces.}
(H. E. Lacey ed.), Univ. of Texas Press, Austin 1980.
%
\bibitem[P.1]{P1} {\bf Pisier, G.}, Holomorphic semi-groups and the 
geometry of Banach spaces. Annals of Math. 115 (1982), 375--392.
%
\bibitem[P.2]{PClev} {\bf Pisier, G.}, On the duality between type and cotype,
{\sl Martingale Theory in Harmonic Analysis and Banach spaces, Proceedings,
Cleveland 1981} (ed. J.A. Chao and W. A. Woyczynski),
Springer Lecture Notes No. 939, 131-144.
%
\bibitem[P.3]{P} {\bf Pisier, G.}, {\sl Volumes of Convex Bodies and 
Banach Spaces Geometry.} Cambridge Univ. Press, 1989.
%
\bibitem[Sz.1]{S1} {\bf Szarek, S. J.}, The finite
dimensional basis
problem with  an appendix on nets of Grassmann manifolds.
{\sl Acta Math.}, 151 (1983), 153--179.
%
\bibitem[Sz.2]{S2} {\bf Szarek, S. J.}, On the existence and
uniqueness of complex structure and spaces with ``few''
operators. {\sl Trans. AMS} 293 (1986), 339--353.
%
\bibitem[T]{T} {\bf Tomczak-Jaegermann, N.}, {\sl
Banach--Mazur Distances
and Finite Dimensional Operator Ideals.} Pitman Monographs and Surveys
in Pure and Applied Mathematics, Longman Scientific \& Technical,
Harlow and John Wiley, New York, 1989.
%
\end{thebibliography}
\end{document}